\documentclass[11pt]{article}
\usepackage{amsfonts,amssymb,amsmath,a4wide,remreset,amscd}

\def\proof {\noindent{\sc{Proof. }}}
\def\qed {\mbox{}\hfill {\small \fbox{}} \\}  
\newcommand{\Zm}{\ensuremath{\mathbb{Z}}}

\newcommand{\Rm}{\ensuremath{\mathbb{R}}}

\newcommand{\Tm}{\ensuremath{\mathbb{T}}}

\newcommand{\Nm}{\ensuremath{\mathbb{N}}}

\def\lto{\longrightarrow}
\def\lmto{\longmapsto}

\def\leq{\leqslant}
\def\geq{\geqslant}

\newcommand{\mI}{\ensuremath{\mathcal{I}}}

\newcommand{\mN}{\ensuremath{\mathcal{N}}}
\newcommand{\mM}{\ensuremath{\mathcal{M}}}

\newcommand{\mA}{\ensuremath{\mathcal{A}}}

\newcommand{\vs}{\vspace{.2cm}}

\makeatletter

\renewcommand{\section}{\@startsection
{section}
{1}
{0mm}
{\baselineskip}
{0.5\baselineskip}
{\normalfont \normalsize \bfseries}}

\renewcommand{\subsection}{\@startsection
{subsection}
{2}
{0mm}
{-\baselineskip}
{0mm}
{\normalfont \normalsize \bfseries}}

\author{Patrick Bernard}

\title{}

\begin{document}
\begin{center}
\begin{scshape}
\begin{Large}
Symplectic aspects of Aubry-Mather theory
\footnote{\today}
\vspace{1cm}\\
\end{Large}
\begin {large}
Patrick Bernard
\footnote{
Institut Fourier, BP 74, 38402, Saint Martin d'H\`eres cedex,
France,\\
patrick.bernard@ujf-grenoble.fr,
http://www-fourier.ujf-grenoble.fr/\~{}pbernard/}
\vspace{1cm}\\\end{large}
\end{scshape}
\end{center}

\begin{small}
\textsc{R\'esum\'e :}
On montre que les ensembles d'Aubry et de Ma\~n\'e introduits
par Mather en dynamique Lagrangienne sont des invariants symplectiques.
On introduit pour ceci une barriere dans l'espace des phases.
Ceci est aussi l'occasion d'\'ebaucher une th\'eorie 
d'Aubry-Mather pour des Hamiltoniens non convexes.
\vspace{1cm}

\textsc{Abstract :}
We prove that the Aubry and Ma\~n\'e sets introduced by Mather
in Lagrangian dynamics are symplectic invariants.
In order to do so, we introduce a barrier on phase space.
This is also an occasion to suggest an Aubry Mather theory 
for non convex Hamiltonians.
\vspace{1cm}\\
\end{small}

In Lagrangian dynamics, John Mather has defined  several invariant sets,
now called the Mather set, the Aubry set, and the Ma\~n\'e set.
These invariant sets provide obstructions to the existence
of orbits wandering in phase space. Conversely, the existence 
of interesting
orbits have been proved 
under some assumptions on the topology of these sets.
Such results were first obtained by John Mather in \cite{MatherFourier},
and then in several papers, see
\cite{Fourier,JAMS,Cheng,Cheng2,Xia1,Xia}
as well as  recent unpublished works of John Mather.

In order to apply these results on examples
one has to understand the topology
of the Aubry and Ma\~n\'e set, which is a very difficult task.
In many perturbative situations, averaging methods appear as
a promising tool in that direction.
In order to use these methods, one has to understand 
how the averaging transformations modify the
Aubry-Mather  sets. 
In the present paper, we answer this question 
and prove  that the Mather set, the Aubry set and the Ma\~n\'e
set are symplectic invariants.

In order to do so, we define a barrier on phase space, which
is some symplectic analogue of the function called the
Peierl's barrier by Mather in \cite{MatherFourier}.
We then propose definitions of Aubry and Ma\~n\'e sets
for general Hamiltonian systems.
We hope that these definitions may also serve as the  starting
point of  an Aubry-Mather theory for some classes of non-convex 
Hamiltonians.
We develop the first steps of such a theory.

Several anterior works gave hints towards the symplectic
nature of Aubry-Mather theory, see \cite{Bulletin, PPS, S1, S2}
for example. 
These works prove the symplectic invariance of the
$\alpha$ function of Mather, and one may consider that 
the symplectic invariance of the Aubry set is not a surprising result
after them. However, the symplectic invariance of the Ma\~n\'e set
is, to my point of view, somewhat unexpected. 
It is possible that the  geometric methods
introduced in \cite{PPS} may also  be used 
to  obtain
symplectic definitions of the Aubry and Ma\~n\'e set.

\section{Mather theory in Lagrangian dynamics}\label{Mather}
We recall the basics of Mather theory
and state our main result, Theorem \ref{main}.
The original references for most of the material presented in this 
section 
are Mather's papers \cite{Mather}
and \cite{MatherFourier}.
The central object is the Peierl's barrier, introduced by Mather
in \cite{MatherFourier}.
Our presentation is also influenced by the work of Fathi \cite{Fathibook}.

\subsection{}\label{HypothesesH}
In this section, we  consider a $C^2$ Hamiltonian function 
$H:T^*M\times \Tm\lto \Rm$, where $M$ is a compact 
connected manifold without boundary, and $\Tm=\Rm/\Zm$.
We  denote by $P=(q,p)$ the points of $T^*M$.
The cotangent bundle is endowed with its canonical one-form
$\eta =pdq$, and with its canonical symplectic form
$\omega=-d\eta$.
Following a very standard device, we  reduce our 
non-autonomous Hamiltonian function $H$
to an autonomous one
by considering  the extended phase space
$T^*(M\times \Tm)=T^*M\times T^*\Tm$.
We denote by $(P,t,E)$, $P\in T^*M$, $(t,E)\in T^*\Tm$
the points of this space.
We consider the canonical one-form $\lambda =pdq+E dt$
and the associated symplectic form  $\Omega=-d\lambda$.
We define the new Hamiltonian
$G:T^*(M\times \Tm)\lto \Rm$
be the expression
$$
G(P,t,E)= E+H(P,t).
$$
We denote by $V_G(P,t,E)$ the Hamiltonian vector-field 
of $G$, which is defined by the relation
$$
\Omega_{(P,t,E)}(V_G,.)=dG_{(P,t,E)}.
$$
We fix once and for all a Riemannian metric
on $M$, and use it to define
norms of tangent vectors and tangent covectors of $M$.
We will denote this norm  indifferently by $|P|$ or by $|p|$
when $P=(q,p)\in T^*_qM$.
We denote by $\pi$ the canonical projections 
$T^*M\lto M$ or $T^*(M\times \Tm)\lto M\times \Tm$.
The theory of Mather relies on the following standard set of hypotheses.
\begin{enumerate}
\item
\textsc{Completeness. } 
The Hamiltonian vector-field $V_G$  on $T^*(M\times \Tm)$ generates
a complete flow, denoted by $\Phi_t$. The flow $\Phi_t$
preserves the level sets of $G$.
\item \textsc{Convexity. } For each $(q,t)\in M\times \Tm$, the
function $p\lmto H(q,p,t)$ is convex on $T_q^*M$,
with positive definite Hessian. Shortly,
$\partial^2_pH >0$.
\item \textsc{Super-linearity. }
For each $(q,t)\in M\times \Tm$, the function $p\lmto H(q,p,t)$
is super-linear, which means that 
$
\lim_{|p|\lto\infty} H(t,x,p)/|p|=\infty.$
\end{enumerate}

\subsection{}\label{HypothesesL}
We associate to the Hamiltonian $H$ a Lagrangian function
$L: TM\times \Tm\lto \Rm$
defined by 
$$L(t,q,v)=\sup _{p\in T^*_q M} p(v)-H(t,q,p).
$$
The Lagrangian satisfies:
\begin{enumerate}
\item \textsc{Convexity. } For each $(q,t)\in M\times \Tm$, the 
function $v\lmto L(q,v,t)$ is a convex function on $T_qM$,
with positive definite Hessian. Shortly, $\partial^2_vL>0$.
\item \textsc{Super-linearity. }
For each $(q,t)\in M\times \Tm$, the function $v\lmto L(q,v,t)$
is super-linear on $T_qM$.
\end{enumerate}
Let $X(t)=(P(t),s+t,E(t))$ be a Hamiltonian orbit of $G$,
and let $q(t)=\pi(P(t))$.
Then we have the identities
$$
\lambda_{X(t)}(\dot X(t))-G(X(t))
=
\eta_{P(t)}(\dot P(t))-H(P(t),s+t)
=
L(q(t), \dot q(t),s+t).
$$
\subsection{}\label{mini}
Following John Mather, we define 
the function $F:M\times \Tm\times M\times \Rm^+ \lto \Rm$ by
$$
F(q_0,t;q_1,s)=
\min _{\gamma} \int _{0}^{s}
L(\gamma(\sigma),\dot \gamma(\sigma),t+\sigma)d\sigma,
$$
where the minimum is taken on the set of absolutely
continuous curves $\gamma:[0,s]\lto  M$ which satisfy
$\gamma(0)=q_0$ and $\gamma(1)=q_1$.
We also define 
the Peierl's barrier 
$h:M\times \Tm \times M\times \Tm\lto \Rm\cup \{\pm \infty\}$ by
$$
h(q_0,t_0;q_1,t_1):= \liminf_{n\in \Nm }  F(q_0,t_0;q_1,s_1+n),
$$
where $t_0+s_1 \text{ mod }1=t_1$.
This barrier is the central object in Mather's study
of globally minimizing orbits.

\subsection{}
Let us set $m(H)=\inf_{(q,t)\in M\times \Tm} h(q,t;q,t)$.
It follows from \cite{Mather}, see also \cite{Mane},
 that $m(H)\in\{-\infty, 0, +\infty\}$.
In addition, for each Hamiltonian $H$
satisfying the hypotheses \ref{HypothesesH}, there exists
one and only one real number $\alpha(H)$ such that
$m(H-\alpha(H))=0$. 
As a consequence, there is no loss of generality in assuming 
that $m(H)=0$, or equivalently that 
$\alpha(H)=0$.
We will make this assumption from now on in this section.
Let us mention the terminology of Ma\~n\'e, who called
super-critical the Hamiltonians $H$ satisfying $m(H)=+\infty$,
sub-critical the Hamiltonians satisfying $m(H)=-\infty$,
and critical the Hamiltonians satisfying $m(H)=0$.

\subsection{}\label{Peierl}
If $m(H)=0$, 
the function $h$ is a real valued Lipschitz function on 
$M\times \Tm \times M\times \Tm$, which
  satisfies the triangle inequality
$$
h(q_0,t_0;q_2,t_2)\leq h(q_0,t_0;q_1,t_1)+h(q_1,t_1;q_2,t_2)
$$
for all $(q_0,t_0)$, $(q_1,t_1)$ and $(q_2,t_2)$ in $M\times \Tm$.
In addition, for each 
$(q,t)\in M\times \Tm$, the function $h(q,t;.,.)$ is a weak KAM solution
in the sense of Fathi, which means that, for $\tau\geq \theta$ in $\Rm$,
and  $x\in M$,
we have 
$$
h(q,t;x,\tau \text{ mod } 1)=\min \Big(
h(q,t;q(\theta),\theta \text{ mod } 1
)+ \int_{\theta}^{\tau} L(q(s),\dot q (s),s)ds \Big)
$$
where the minimum is taken on the set of absolutely continuous
curves $q(s):[\theta,\tau]\lto M$ such that $q(\tau)=x$. 
Similarly, we have,  for $\tau\geq \theta$ in $\Rm$,
and  $x\in M$,
$$
h(x,\theta \text{ mod } 1;q,t)=\min \Big(
h(q(\tau),\tau \text{ mod } 1;q,t)+ 
\int_{\theta}^{\tau} L(q(s),\dot q (s),s)ds \Big)
$$
where the minimum  is taken on the set of absolutely continuous
curves $q(s):[\theta,\tau]\lto M$  such that $q(\theta )=x$.

\subsection{}\label{Aubry}
The projected Aubry set $\mA(H)$ is the set of points 
$(q,t)\in M\times \Tm$ such that 
$h(q,t;q,t)=0$.
Albert Fathi proved that, for each point $(q,t)\in \mA(H)$, the function
$h(q,t;.,.)$ is differentiable at $(q,t)$.
Let us denote by $X(q,t)$
the differential $\partial_3h(q,t;q,t)\in T_q^*M$
of the function
$h(q,t;.,t)$ at point $q$.
The Aubry set $\tilde \mA(H)$
is defined as
$$
\tilde \mA(H)=\{(X(q,t),t,-H(X(q,t),t));  
(q,t)\in \mA(H)\}\subset T^*(M\times \Tm).
 $$
The Aubry set is compact, $\Phi$-invariant, and it is a Lipschitz
graph over the projected Aubry set $\mA(H)$.
These are results of John Mather, see \cite{MatherFourier}.
In our presentation, which follows Fathi, this amounts
to say that the function $(q,t)\lmto X(q,t)$ is Lipschitz on
$\mA(H)$.

\subsection{}
The Mather set $\tilde \mM(H)$ is defined as the union of the supports
of all $\Phi$-invariant probability
measures on $T^*(M\times \Tm)$
 concentrated on $\tilde \mA(H)$.
This set was first defined by Mather, but our definition
is due to Ma\~n\'e.

\subsection{}\label{pMane}
The projected Ma\~n\'e set $\mN(H)$ is the set of points
 $(q,t)\in M\times \Tm$
such that there exist points $(q_0,t_0)$ and $(q_1,t_1)$ in $\mA(H)$,
satisfying
$$
h(q_0,t_0;q_1,t_1)=h(q_0,t_0;q,t)+h(q,t;q_1,t_1).
$$
Let us denote by $\mI(q_0,t_0;q_1,t_1)$
the set of points $(q,t)\in M\times \Tm$ which satisfy this relation.
If $(q_0,t_0)\in \mA(H)$ 
and $(q_1,t_1) \in \mA(H)$ are given, and if 
$(q,t)\in \mI(q_0,t_0;q_1,t_1)$,
then
the function $h(q_0,t_0;.,t)$ is differentiable at $q$, 
as well as the function $h(.,t;q_1,t_1)$, and 
$\partial_3h(q_0,t_0,q,t)+\partial_1h(q,t;q_1,t_1)=0$. 
This is proved in \cite{JAMS}
following ideas of Albert Fathi.
We define
$$
\tilde \mI (q_0,t_0;q_1,t_1)
:=\Big\{\Big(\partial_3h(q_0,t_0,q,t), t,-H\big(\partial_3h(q_0,t_0,q,t),t\big)
\Big)
, (q,t)\in  \mI (q_0,t_0;q_1,t_1)\Big\}.
$$ 
The set $\tilde \mI (q_0,t_0;q_1,t_1)$ is a compact 
$\Phi$-invariant subset of $T^*(M\times \Tm)$, and it is 
a Lipschitz Graph. 
The Ma\~n\'e set $\tilde \mN(H)$
is the set 
$$
\tilde \mN(H)
=\bigcup_{(q_0,t_0),(q_1,t_1)\in \mA(H)}
\tilde \mI(q_0,t_0;q_1,t_1)
\subset T^*(M\times \Tm).
$$
The Ma\~n\'e set was first introduced by Mather in \cite{MatherFourier},
it is compact and $\Phi$-invariant, and  it contains the Aubry set.
In other words, we have the important inclusions
$$
\tilde \mM(H)\subset \tilde \mA(H)
\subset \tilde \mN(H).
$$
The Ma\~n\'e set is usually not a graph. However,
it satisfies 
$$
\tilde \mN(H)\cap \pi^{-1}\big( \mA(H)\big)
=\tilde \mA(H).
$$
This follows from the fact, proved by Albert Fathi,
that, for each $(x,\theta)\in M\times \Tm$ and
each $(q,t)\in \mA(H)$, the function
$h(x,\theta;.,t)$ is differentiable at $q$
and satisfies 
$\partial_3h(x,\theta;q,t)=X(q,t)$.

\subsection{}
Mather introduced the function 
$d(q,t;q',t')=h(q,t;q',t')+h(q',t';q,t)$ on $M\times \Tm$.
When restricted to $\mA(H) \times \mA(H)$, it  is
a pseudo-metric. This means that this function is symmetric, 
non-negative,
satisfies the triangle inequality, and $d(q,t;q,t)=0$ for 
$
(q,t)\in \mA(H)$.
We shall also denote by $d$ the pseudo-metric
$d(P,t,-H(P,t);P',t',-H(P',t'))=d(\pi(P),t;\pi(P'),t')$ on $\tilde A(H)$.
The relation $d(P,t,E;P',t',E')=0$ is an equivalence relation 
on $\tilde \mA(H)$.
The classes of equivalence are called the static classes.
Let us denote by $\dot \mA(H)$ the set of static classes.
The pseudo-metric $d$ gives rise to a metric
$\dot d$ on  $\dot \mA(H)$.
The compact metric space 
$(\dot \mA(H),\dot d)$ is called the quotient Aubry set.
It was introduced by John Mather.

\subsection{}\label{main}
The diffeomorphism $\Psi:T^*(M\times \Tm)\lto T^*(M\times \Tm)$
is called exact if the form  $\Psi^*\lambda-\lambda$
is exact.
\vs

\noindent\textsc{Theorem}
\begin{itshape}
Let $H$ be a Hamiltonian satisfying the hypotheses \ref{HypothesesH},
and let  $\Psi:T^*(M\times\Tm)\lto 
T^*(M\times \Tm)$ be an exact diffeomorphism
such that the Hamiltonian 
$$
\Psi^*H:= G\circ\Psi (P,t,E)-E
$$
is independent of $E$ and satisfies the hypotheses \ref{HypothesesH}
when considered as a function on $T^*M\times \Tm$.
Then $m(\Psi^*H)=m(H)$ hence 
$\alpha(H)=\alpha(\Psi^*H)$. If $m(H)=0$, then
we have 
$$
\Psi(\tilde \mM(\Psi^*H))=\tilde \mM(H)\,,\;\;
\Psi(\tilde \mA(\Psi^*H))=\tilde \mA(H)\,,\;\;
\Psi(\tilde \mN(\Psi^*H))=\tilde \mN(H).
$$
In addition, $\Psi$ sends the static classes of $\Psi^*H$
onto the static classes of $H$, and the induced mapping
$$\dot \Psi:\dot \mA(\Psi^*H)\lto \dot \mA(H)$$
is an isometry for the quotient metrics.
\end{itshape}

\subsection{}
We prove this result in the sequel.
In section \ref{barrier}, we set the basis of a symplectic Aubry-Mather
theory for general Hamiltonian systems. 
We prove that the analogue of Theorem \ref{main}
holds in this general setting.
We also continue the theory a bit further
than would be necessary to prove
Theorem \ref{main}.
In section \ref{fin}, we prove that, under the hypotheses
of Theorem \ref{main},
the symplectic Aubry-Mather sets coincide with the standard
Aubry-Mather sets, which ends the proof of Theorem \ref{main}.

\section{A barrier in phase space}\label{barrier}
%
%
%
%
We propose general definitions for a Mather theory
of Hamiltonian systems. Of course, the definitions
given below provide relevant objects only for some
specific Hamiltonian systems. It would certainly be interesting to
give natural conditions on $H$ implying non-triviality 
of the theory developed in this section. We shall only
check, in the next section, that our definitions
coincide with the standard ones in the convex case, obtaining
non-triviality in this special case.
Let us mention once again that it might be possible and interesting
to find more geometric definition using the 
methods of \cite{PPS}.

\subsection{}
In this section, we work in a very general setting.
We consider a manifold $N$, not necessarily compact,
and an autonomous Hamiltonian function 
$G:T^*N\lto \Rm$. 
We assume that $G$ generates
a complete Hamiltonian flow $\Phi_t$.
We make no convexity assumption.
We denote by $\lambda$ the canonical one-form of $T^*N$,
and by $V_G(P)$ the Hamiltonian vector-field of $G$.
Let $D(P,P')$ be a distance on $T^*N$ induced from a Riemannian 
metric. We identify $N$ with the zero section of $T^*N$,
so that $D$ is also a distance on $N$. We assume that
$D(\pi(X),\pi(X'))\leq D(X,X')$ for $X$ and $X'$ in $T^*N$.

  \subsection{}
Let $X_0$ and $X_1$ be two points of $T^*N$.
A pre-orbit between $X_0$ and $X_1$ is the data of a sequence
$\underline Y=(Y_n)$ of curves $Y_n(s):[0,T_n]\lto T^*N$ such that:
\begin{enumerate}
\item
For each $n$, 
the curve $Y_n$ has a finite number $N_n$ of discontinuity
points  $T^i_n\in]0,T_n[,1\leq i\leq N_n$ such that 
$T^{i+1}_n>T^i_n$. We shall also often use the notations
$T^0_n=0$ and $T^{N_n+1}_n=T_n$.
\item
The curve $Y_n$ satisfies 
$Y_n(T^i_n+s)=\Phi_s(Y_n(T^i_n))$ for each 
$s\in [0,T^{i+1}_n-T^i_n[$.
We denote by $Y_n(T^i_n-)$ the point 
$\Phi_{T^i_n-T^{i-1}_n}(Y(T^{i-1}_n))$ and impose that 
$Y_n(T_n)=Y_n(T_n-)$.
\item
We have  $T_n\lto \infty$ as $n\lto \infty$.
\item
We have  $Y_n(0)\lto X_0$ and $Y_n(T_n)\lto X_1$. In addition, we have 
$\lim_{n\lto \infty}\Delta(Y_n)=0$, where 
we denote by $\Delta(Y_n)$ the sum 
$\sum_{i=1}^{N_n} D(Y_n(T^i_n-),Y_n(T^i_n))$.
\item
There exists a compact subset $K\subset T^*N$ which contains 
the images of all the curves $Y_n$.
\end{enumerate}
The pre-orbits do not depend on the metric which has been used
to define the distance $D$.
In a standard way, we call action of the curve $Y_n(t)$
the value 
$$
A(Y_n)=\int_0^{T_n}\lambda _{Y_n(t)}(\dot Y_n(t))-G(Y_n(t))\, dt.
$$
The action of the pre-orbit $\underline Y$ is 
$$
A(\underline Y):= \liminf_{n\lto \infty} A(Y_n).
$$

\subsection{}
\textsc{Lemma}
\begin{itshape}
If there exists a pre-orbit between $X_0$ and $X_1$,
then $G(X_0)=G(X_1)$.\vs\\
\end{itshape}
\proof
This follows easily from the fact that the Hamiltonian flow
$\Phi$ preserves the Hamiltonian function $G$.
\qed
\subsection{}
We define the barrier $\tilde h:T^*M \times T^*M \lto \Rm \cup \{\pm \infty\}$
by the expression
$$
\tilde h(X_0,X_1)=\inf _{\underline Y} A(\underline Y)
$$
where the infimum is taken on the set of  pre-orbits between $X_0$ 
and $X_1$.
As usual, we set $  \tilde h(X_0,X_1)=+\infty$ if there does not exist
any pre-orbit between $X_0$ and $X_1$.
If $\tilde h(X_0,X_1)<+\infty$, then the forward orbit of 
$X_0$ and the backward orbit of $X_1$ are bounded.
As a consequence, if $\tilde h(X,X)<+\infty$, then the orbit of $X$
is bounded.

\subsection{}\label{invariance}                
\textsc{Property}
\begin{itshape}
For each $t>0$, we  have the equality
$$
\tilde h(X_0,X_1)=
\tilde h(\Phi_t(X_0),X_1)
+\int_0^t 
\lambda_{\Phi_s(X_0)}
\big(V_G(\Phi_s(X_0))
\big)-
G(\Phi_s(X_0))\,ds
$$
and 
$$
\tilde h(X_0,\Phi_t(X_1))=
\tilde h(X_0,X_1)
+\int_0^t 
\lambda_{\Phi_s(X_1)}
\big(V_G(\Phi_s(X_1))
\big)-
G(\Phi_s(X_1))\,ds
$$

\end{itshape}
\proof
We shall prove the first equality, the proof of the second one is similar.
To each pre-orbit $\underline Y$ between $X_0$
and $X_1$, we associate the pre-orbit $\underline Z$
between $\Phi_t(X_0)$ and $X_1$
defined by $Z_n(s):[0,T_n-t]\lmto Y_n(s+t)$.
We have 
$$
A(\underline Y)=A(\underline Z)+
\int_0^t
\lambda_{\Phi_s(X_0)}
\big(V_G(\Phi_s(X_0))
\big)-
G(\Phi_s(X_0))\,ds
$$
This implies that
$$
\tilde h(\Phi_t(X_0),X_1)\leq
\tilde h(X_0,X_1)-
\int_0^t \lambda_{\Phi_s(X_0)}
\big(V_G(\Phi_s(X_0))
\big)-
G(\Phi_s(X_0))\,ds.
$$
In a similar way, we associate to each pre-orbit 
$\underline Z=Z_n(s):[0,T_n]\lto T^*M$
between $\Phi_t(X_0)$ and $X_1$  the  pre-orbits
$\underline Y:[0,T_n+t]\lto T^*M$ between $X_0$ and $X_1$
defined by $Y_n(s)=\Phi_{s-t}(Z_n(0))$ for $s\in [0,t]$
and
$Y_n(s)=Z_n(s-t)$ for $s\in [t,T_n+t]$.
We have 
$$A(\underline Y)=A(\underline Z)
+\int_0^t \lambda_{\Phi_s(X_0)}
\big(V_G(\Phi_s(X_0))
\big)-
G(\Phi_s(X_0))\,ds.
$$
This implies that 
$$
\tilde h(X_0,X_1)\leq
\tilde h(\Phi_t(X_0),X_1)+
\int_0^t \lambda_{\Phi_s(X_0)}
\big(V_G(\Phi_s(X_0))
\big)-
G(\Phi_s(X_0))\,ds.
$$
\qed

\subsection{}
\textsc{property}
\begin{itshape}
The function $\tilde h$ satisfies the triangle inequality. 
More precisely, the relation
$$
\tilde h(X_1,X_3)\leq \tilde    h(X_1,X_2)+  \tilde    h(X_2,X_3)
$$
holds for each points $X_1$, $X_2$ and $X_3$ such that the right hand side
has a meaning. \vs\\
\end{itshape}
\proof
If one of the values $\tilde    h(X_1,X_2)$ or  $\tilde    h(X_2,X_3)$
is $+\infty$, then there is nothing to prove.
If they are both different from $+\infty$, then,
for each $\epsilon>0$ there exists
a  pre-orbits $\underline Y=Y_n:[0,T_n]\lto T^*N$
 between $X_1$ and $X_2$ such that
$A(\underline Y)\leq \tilde    h(X_1,X_2)+\epsilon$ 
(resp. 
$A(\underline Y)\leq -1/\epsilon$ 
in the case where $ \tilde    h(X_1,X_2)=-\infty$)
and 
a  pre-orbits $\underline Y'=Y_n':[0,S_n]\lto T^*N$
 between $X_2$ and $X_3$ such that
$
A(\underline Y')\leq \tilde    h(X_2,X_3)+\epsilon$
(resp. 
$A(\underline Y')\leq -1/\epsilon$ 
in the case where $ \tilde    h(X_1,X_2)=-\infty$).
Let us consider the sequence of curves $Z_n(t):[0,T_n+S_n]\lto T^*N$
such that $Z_n=X_n$ on $[0,T_n[$ and 
$Z_n(t+T_n)=Y_n(t)$ for $t\in [0,S_n]$.
It is clear that the sequence $\underline Z=Z_n$ is a  pre-orbit
between $X_1$ and $X_3$, and that its action satisfies 
$$
A(\underline Z)
=A(\underline X)+A(\underline Y)\leq 
\tilde    h(X_1,X_2)+  \tilde    h(X_2,X_3)+2\epsilon.
$$
As a consequence, for all $\epsilon >0$, we have
$\tilde h(X_1,X_3)\leq 
\tilde    h(X_1,X_2)+  \tilde    h(X_2,X_3)+2\epsilon$
hence  the triangle inequality holds. 
\qed

\subsection{}\label{symplecticinvariance}
\textsc{property}
\begin{itshape}
Let $\Psi:T^*N\lto T^*N$ be an exact diffeomorphism.
We have the equality
$$
\tilde h_{G\circ \Psi}(X_0,X_1)=
\tilde h_G(\Psi(X_0),\Psi(X_1))+S(X_0)-S(X_1),
$$
where $S:T^*N\lto \Rm$ is a function such that
$\Psi^*\lambda-\lambda=dS$.\vs
\end{itshape}

\proof
Observe first that $\underline Y=Y_n$
is a pre-orbit for the Hamiltonian $G\circ \Psi$
between points $X_0$ and $X_1$ if and only if 
$\Psi(\underline Y)=\Psi(Y_n)$ is a  pre-orbit
for the Hamiltonian $G$ between $\Psi(X_0)$ and
$\Psi(X_1)$.
As a consequence, it is enough to prove that 
$$
A_{G\circ \Psi}(\underline Y)=A_G(\Psi(\underline Y))
+S(X_0)-S(X_1).
$$
Let us denote by $\underline Z=Z_n$ the pre-orbit $\Psi(Y_n)$.
Setting $T^0_n=0$ and $T^{N_n+1}_n=T_n$, we have  
$$
A_G(Z_n)=
\sum _{i=0}^{N_n}
\int_{T^i_n}^{T^{i+1}_n}
\lambda_{Z_n(t)}(\dot Z_n(t))
-G(Z_n(t))
dt
$$
$$
=\sum _{i=0}^{N_n}
\int_{T^i_n}^{T^{i+1}_n}
(\Psi^*\lambda)_{Y_n(t)}( \dot Y_n(t))
-G\circ \Psi(Y_n(t))
dt
$$
$$
=\sum _{i=0}^{N_n}\left(
\int_{T^i_n}^{T^{i+1}_n}
\lambda _{Y_n(t)}( \dot Y_n(t)) 
-G\circ \Psi(Y_n(t))
dt
+S(Y_n(T^{i+1}_n-))-S(Y_n(T^i_n))
\right)
$$
$$
=A_{G\circ \Psi}(Y_n)
-S(Y_n(0))+
S(Y_n(T_n))+
\sum _{i=1}^{N_n}
\big(
S(Y_n(T^{i}_n-))-S(Y_n(T^i_n)).
\big)
$$
Since the function $S$ is Lipschitz on the compact set $K$ which contains
the image of the curves $Y_n$,
we obtain at the limit
$$A_G(\underline Z)=
A_{G\circ\Psi}(\underline Y)-S(X_0)+
S(X_1).
$$
\qed

\subsection{}
\textsc{proposition}
\begin{itshape}
Let us set $\tilde m(H):=\inf_{X\in T^*N} \tilde h(X,X)$.
We have $\tilde m(H)\in \{-\infty,0,+\infty\}$.
In addition, if $\tilde m(H)=0$, then there exists a point $X$ in $T^*N$
such that $\tilde h(X,X)=0$.\vs
\end{itshape}

\proof
It follows from the triangle inequality that, for each $X\in T^*N$, 
$\tilde h(X,X)\geq 0$ or $\tilde h(X,X)=-\infty$.
As a consequence, $\tilde m(H)\geq 0$ or $\tilde m(H)=-\infty$.
Let us assume that $\tilde m(H)\in [0,\infty[$.
Then there exists a point $X_0\in T^*N$
and a  pre-orbits
$\underline Y=Y_n:[0,T_n]\lto T^*N$ between $X_0$ and $X_0$ such that 
$A(\underline Y)\in [0,\infty[$.
Let $K$ be a compact subset of $T^*N$ which contains 
the image of all the curves $Y_n$.
Let $S_n$ be a sequence of integers such that 
$T_n/S_n\lto \infty$ and  $S_n\lto \infty$.
Let $b_n$ be the integer part of $T_n/S_n$.
Note that $b_n\lto \infty$.
Let $d_n$ be a sequence of integers such that 
$d_n\lto \infty$ and $d_n/b_n\lto 0$.
Since the set $K$ is compact, there exists a sequence
$\epsilon_n\lto 0$ such that, whenether $b_n$ points are given
in $K$, then at least $d_n$ of them lie in a same ball
of radius $\epsilon_n$.
So there exists a point $X_n\in K$ such that at least $d_n$
of the points $Y_n(S_n), Y_n(2S_n), \ldots, Y_n(b_nS_n)$
lie in the ball of radius $\epsilon_n$ 
and center $X_n$.
Let us denote by $Y_n(t_n^1),Y_n(t_n^2) \ldots,Y_n( t_n^{d_n})$
these points, where $t^{i+1}_n\geq t^i_n+S_n$.
We can assume, taking a subsequence, that the sequence
$X_n$ has a limit $X$ in $K$.
It is not hard to see that 
$\underline Y^i=Y_{n|[t^i_n,t^{i+1}_n]}$ is a 
pre-orbit between $X$ and $X$.
On the other hand, for each  $k\in \Nm$, 
we define the sequence of curves
 $Z^k_n:[0,T_n+t^1_n-t^{k}_n]\lto T^*N$
by $Z^k_n(t)=Y_n(t)$ for $t\in[0,t^1_n[$,
and $Z^k_n(t)=Y_n(t+t^{k}_n-t^1_n)$ for $t\in[t^1_n,T_n+t^1_n-t^{k}_n]$.
For each $k$, the  sequence $Z^k_n$ is a  pre-orbit between $X_0$ and $X_0$.
We have 
$$
A(Y_n)=A(Z^k_n)+\sum_{i=1}^{k-1}A(Y^i_n)
$$
hence 
$$
A(\underline Y)\geq \tilde h(X_0,X_0)+(k-1)\tilde h(X,X).
$$
Since $A(\underline Y)$ is a real number,
and since this inequality holds for all $k\in \Nm$, this implies that 
$\tilde h(X,X)=0$.
\qed

\subsection{}
Let us define the symplectic Aubry set of $G$ as the set 
$$
\tilde \mA_s(G):=
\{X\in T^*N \text{ such that } \tilde h(X,X)=0 \text{ and } G(X)=0\}
\subset T^*N.
$$
The symplectic Mather set $\tilde \mM_s(G)$ of $G$
is the union of the supports of the compactly supported 
$\Phi$-invariant probability measures concentrated on $\tilde \mA_s(G)$.
Note that, in general, it is not clear that the symplectic
Aubry set should be closed. The symplectic Mather set,
then, may not be contained in the symplectic Aubry set,
but only in its closure.
The Mather set and the Aubry  set are $\Phi$-invariant,
as follows directly from \ref{invariance}.
If $\tilde m(H)=0$, then the symplectic  Aubry set is not empty,
and all its orbits are bounded, 
hence the symplectic Mather set $\tilde \mM_s(G)$
is not empty.

\subsection{}
For each pair $X_0$, $X_1$ of points in $\tilde \mA_s(G)$,
we define the set 
$\tilde \mI_s(X_0,X_1)$ of points $P\in T^*N$
such that 
$$
\tilde h(X_0,X_1)=\tilde h(X_0,X)+\tilde h(X,X_1)
$$
if $\tilde h(X_0,X_1)\in \Rm$, and 
$\tilde \mI_s(X_0,X_1)=\emptyset$ otherwise.
Note that the sets $\tilde  \mI_s(X_0,X_1)$
are all contained in the level $\{G=0\}$.
Indeed, the finiteness of $\tilde h(X_0,X)$
implies that $G(X_0)=G(X)$, while $G(X_0)=0$
by definition of $\tilde \mA_s(G)$.
It follows from \ref{invariance} that the set 
$\tilde \mI_s(X_0,X_1)$ is $\Phi$-invariant.
We now define the symplectic Ma\~n\'e set as
$$
\tilde \mN_s(G):=
\bigcup_{X_0,X_1\in \tilde \mA_s(G)}
\tilde \mI_s(X_0,X_1).
$$
The  Ma\~n\'e set is $\Phi$-invariant,
all its orbits are bounded.
We have the inclusion 
$$
\tilde \mA_s(G)\subset \tilde \mN_s(G).
$$
In order to prove this inclusion, just observe that 
$
X_0\in \tilde \mI(X_0,X_0)
$
for each $X_0\in \tilde \mA_s(G)$.

\subsection{}
If $\Psi:T^*N\lto T^*N$ is an exact diffeomorphism,
then we have 
$$
\Psi(\tilde  \mM_s(G\circ \Psi))=\tilde \mM_s(G),\;\;
\Psi(\tilde \mA_s(G\circ \Psi))=\tilde \mA_s(G),\;\;
\Psi(\tilde \mN_s(G\circ \Psi))=\tilde \mN_s(G),
$$
this follows obviously from \ref{symplecticinvariance},
and from the fact that $\Psi$ conjugates the Hamiltonian flow
of $G$ and the Hamiltonian flow of $G\circ \Psi$.

\subsection{}
Let us assume that $\tilde m(G)=0$, and set
$$
\tilde d(X,X')=\tilde h(X,X')+\tilde h(X',X).
$$
We have $\tilde d(X,X')\geq 0$, and the function $\tilde d$
satisfies the triangle inequality, and is symmetric.
In addition, we obviously have 
$\tilde d(X,X)=0$ if and only if $X\in \tilde \mA_s(G)$.
The restriction of the function $\tilde d$ to the set 
$\tilde \mA_s(G)$ is a pseudo-metric with $+\infty$
as a possible value.
We define an equivalence relation on $\tilde \mA_s(G)$
by saying that the points $X$ and $X'$ are equivalent if and only if
$\tilde d(X,X')=0$.
The equivalence classes of this relation are called the static classes.
Let us denote by $(\dot \mA_s(G),\dot d_s$)
the metric space obtained from $\tilde \mA_s$ by identifying
points $X$ and $X'$ when $\tilde d(X,X')=0$.
In other words, the set $\dot \mA_s(G)$ is the set of static classes of $H$.
We call $(\dot \mA_s(G),\dot d_s)$
the quotient Aubry set. Note that the metric 
$\dot d_s$ can take the value $+\infty$.
The quotient Aubry set is also well behaved under
exact diffeomorphisms.
More precisely, if $\Psi$ is an exact diffeomorphism of $T^*N$,
then the image of a static class of $G \circ \Psi$
is a static class of $G$.
This defines a map
$$
\dot \Psi :\dot \mA_s(G\circ\Psi )
\lto \dot\mA_s(G)
$$
which is an isometry for the quotient metrics.

\subsection{}
\textsc{Proposition, }
\begin{itshape}
Assume that  $\tilde m(G)=0$, and in addition
that the function $\tilde h$ is bounded from below.
Then the orbits of $\tilde \mN_s(G)$
are bi-asymptotic to $\tilde \mA_s(G)$.
In addition, 
for each orbit $X(s)$ in $\tilde \mN_s(G)$,
there exists a static class $S-$ in $\tilde \mA_s(G)$
and a static class $S+$  such that the orbit $X(s)$
is $\alpha$-asymptotic to $S-$ and $\omega$-asymptotic to $S+$.\vs
\end{itshape}

\proof
Let $\omega$ and $\omega'$ be two points in the $\omega$-limit
of the orbits $X(t)=\Phi_t(X)$.
We have to prove that $\omega$ and $\omega'$ belong to the
symplectic Aubry set, and to the same
static class.
It is enough to prove that 
$\tilde d(\omega,\omega')=0$.
In order to do so, we consider two increasing sequences 
$t_n$ and $s_n$, such that 
$t_n-s_n\lto \infty$, $s_n-t_{n-1}\lto \infty$,
$X(t_n)\lto \omega$ and $X(s_n)\lto \omega'$.
Let $\underline Y=Y_n:[0,t_n-s_n]\lto T^*N$
be the pre-orbit between $\omega'$ and $\omega$ 
defined by
$Y_n(t)=X(t-s_n)$.
Similarly, we consider the pre-orbit 
$\underline Z =Z_n:[0,s_{n+1}-t_n]\lto T^*N$
between $\omega$ and $\omega'$
defined by 
$Z_n(t)=X(t-t_n)$.
Since $X$ belongs to $\tilde \mN_s(G)$, there exist
points $X_0$ and $X_1$ in $\tilde \mA_s(G)$ such that
$X\in \tilde \mI(X_0,X_1)$.
In view of \ref{invariance}, we have
$$
\tilde h(X(t_n),X_1)=\tilde h(X(t_m),X_1)
+\int_{t_n}^{t_m} \lambda_{X(t)} (\dot X(t))
-G(X(t))dt
$$
for all $m\geq n$. 
Since the function
$\tilde h$
is bounded from below, we conclude that 
 the double sequence 
$\int_{t_n}^{t_m} \lambda_{X(t)} (\dot X(t))
-G(X(t))dt, m\geq n
$
is bounded from above, so that 
$$
\liminf
\int_{t_n}^{t_{n+1}} \lambda_{X(t)} (\dot X(t))
-G(X(t))dt
\leq 0.
$$
As a consequence, we have 
$
\liminf A(Y_{n+1}) +A(Z_n)
\leq 0
$
hence $A(\underline Y)+A(\underline Z)=0$,
and $\tilde d(\omega,\omega')=0$.
The proof is similar for the $\alpha$-limit.
\qed
It is useful to finish with section with a technical remark.

\subsection{}\label{tech}
\textsc{Lemma}
\begin{itshape}
Let $\underline Y=Y_n:[0,T_n]\lto T^*N$
be a pre-orbit between 
between $X_0$ and $X_1$.
There exists a pre-orbit $\underline Z$
between  $X_0$ and $X_1$ which has the same action as $\underline Y$,
and has discontinuities only at times
$1,2, \ldots, [T_n]-1$, where $[T_n]$ is the integer part
of $T_n$.\vs
\end{itshape}

\proof
We set 
$Z_n(k+s)=\Phi_s(Y_n(k))$ for each $k= 0, 1,\ldots, [T_n]-2$,
and $s\in [0,1[$, and $Z_n([T_n]-1+s)=\Phi_s(Y_n([T_n]-1))$
for each $s\in [0,1+T_n-[T_n][$.
It is not hard to see that $A(Z_n)-A(Y_n)\lto 0$,
hence $A(\underline Y)=A(\underline Z)$.
\qed
%
%
%
%
%
\section{The case of convex Hamiltonian systems}\label{fin}
%
%
%
%
We assume the hypotheses 
\ref{HypothesesH}, and prove that the symplectic definitions
of section \ref{barrier}   agree with the standard definitions
of section \ref{Mather}.
This proves that the theory of section \ref{barrier} is not trivial
at least in this case. This also ends the proof of Theorem \ref{main}.

\subsection{}
In this section, we consider a Hamiltonian function
$H:T^*M\times \Tm\lto\Rm$ satisfying the hypotheses \ref{HypothesesH}.
We set $N=M\times \Tm$.
We denote by $(P,t,E)$ the points of 
$T^*N$ and  set $G(P,t,E)=E+H(P,t):T^*N\lto \Rm$.
We denote by $h(q,t;q',t')$ the Peierl's barrier associated to 
$H$ in section \ref{Mather}
and by $\tilde h(P,t,E;P',t',E')$ the barrier associated to $G$ in section 
\ref{barrier}.

\subsection{}
Before we state the main result of this section, some
terminology is necessary.
If $u:M\lto \Rm$ is a continuous function, we say that
$P\in T_q^*M$ is a proximal super-differential
of $u$ at point $q$
(or simply a super-differential)
if there exists a smooth function $f:M\lto \Rm$ 
such that $f-u$ has a minimum at $q$ and $df_q=P$.
Clearly, if $u$ is differentiable at $q$ and if
$P$ is a proximal super-differential of $u$ at $q$,
then $P=du_q$.

\subsection{}
\textsc{Proposition}
\begin{itshape}
We have the relation
$$
h(q,t;q',t')=
\min_{P\in T_q^*M,
P'\in T_{q'}^*M} \tilde h(P,t,-H(P,t);P',t',-H(P',t')).
$$
In addition, if  the minimum is reached at  $(P ,P')$
then $P$ is a super-differential of the function 
$h(.,t;q',t')$
at point $q$
and  $-P'$ is a super-differential of the function 
$h(q,t;.,t')$ at point $q'$.\vs\\
\end{itshape}
\proof
Let us fix two points $(q,t)$ and $(q',t')$
in $N=M\times \Tm$.
We claim that the inequality
$$
\tilde h(P,t,E;P',t',E')\geq h(q,t;q',t')
$$
holds for each $(P,t,E)\in T_{(q,t)}^*N$
and each  $(P',t',E')\in T_{(q',t')}^*N$.
If $\tilde h(P,t,E;P',t',E')=+\infty$, 
then there is nothing to prove.
Else, let us fix $\epsilon>0$.
There exists a pre-orbit $\underline Y=Y_n(s):[0,T_n]\lto T^*N$
between $(P,t,E)$ and $(P',t',E')$ 
such that $A(\underline Y)\leq \tilde h(P,t,E;P',t',E')+\epsilon$
(resp.  $A(\underline Y)\leq -1/\epsilon$
in the case where $ \tilde h(P,t,E;P',t',E')=-\infty$).
In view of \ref{tech}, it is possible to assume that
the discontinuity 
points $T_n^i$ of $Y_n$ satisfy $T_n^{i+1}\geq T_n^i+1$.
Let us write 
$$
Y_n(s)=(P_n(s),\tau_n(s),E_n(s)),
$$
and $q_n(s)=\pi(P_n(s))$.
Let $\delta_n^i$ be the real number closest to 
$T_n^{i+1}-T^i_n$ among those which satisfy
$\tau_n(T_n^i)+\delta_n^i=\tau_n(T_n^{i+1})$.

We have 
$$
A(Y_n)=
\sum_{i=0}^{N_n}
\int_{T_n^i}^{T_n^{i+1}}
L(q_n(s),\dot q_n(s),s+\tau_n(T_n^i))dt
\geq \sum_{i=0}^{N_n}
F(q(T^i_n),\tau_n(T^i_n);q(T^{i+1}_n-),T_n^{i+1}-T_n^i).
$$
It is known that the functions $F(q,t;q',s)$
is Lipschitz on $\{s \geq 1\}$, see for example \cite{Fourier}, 3.2.
We have
$$
\sum _{i=0}^{N_n} \Big|
F(q_n(T^i_n),\tau_n(T^i_n);q_n(T^{i+1}_n-),T_n^{i+1}-T_n^i)
-
F(q_n(T^i_n),\tau_n(T_n^i);q_n(T^{i+1}_n),\delta_n^i)
\Big|
$$
$$
\leq C\sum _{i=0}^{N_n-1}D(q_n(T^{i+1}_n-),\tau_n(T^{i+1}_n-);
q_n(T^{i+1}_n),\tau_n(T^{i+1}_n))
$$
$$
\leq C\sum _{i=0}^{N_n-1}D(Y_n(T^{i+1}_n-),Y_n(T^{i+1}_n))
\lto 0.
$$
As a consequence, we have 
$$
A(\underline Y)\geq 
\liminf \sum_{i=0}^{N_n}
F(q(T^i_n),\tau_n(T^i_n);q(T^{i+1}_n),\delta_n^i)
$$
$$
\geq \liminf F\left(q_n(0),\tau_n(0);q_n(T_n), \sum_{i=0}^{N_n}\delta_n^i\right)
\geq h(q,t;q',t'),
$$
hence 
$\epsilon +\tilde h(P,t,E;P',t',E')\geq h(q,t;q',t')$
(resp. $-1/\epsilon\geq  h(q,t;q',t')$).
Since this holds for all $\epsilon>0$, we have 
$\tilde h(P,t,E;P',t',E')\geq h(q,t;q',t')
$
as desired.

Conversely, let us consider
a sequence $T_n$ such that $T_n\lto \infty$,
$t+T_n\text{ mod } 1=t'$, and  
$$
h(q,t;q',t')=\lim_{n\lto \infty} F(q,t;q',T_n).
$$
Let $q_n(s):[0,T_n]\lto M$ be a curve
such that
$$
\int_0^{T_n}
L(q_n(s),\dot q_n(s), s+t) ds
= F(q,t;q',T_n).
$$
Since the curve $q_n$ is minimizing the action,
there exists a Hamiltonian trajectory  
$$
Y_n(s)=(P_n(s), t +s,E_n(s)=-H(X_n(s),t+s))
:[0, T_n]\lto T^*N
$$
 whose projection on 
$M$ is the curve $q_n$.
In addition,  by well known results on minimizing orbits, 
see \cite{Mather}, 
there exists a compact subset of $T^*M$ which contains
the images of all the curves  $P_n(s)$.
As a consequence, we can assume, taking a subsequence if necessary, 
that the sequences $P_n(0)$ and $P_n(T_n)$ have limits
$P\in T_q^*M$ and $P'\in T_{q'}^*M$.
The sequence $\underline Y=Y_n$
is then a pre-orbit between $(P,t,-H(P,t))$ and 
$(P',t',-H(P',t'))$,
and its action is 
$$
A(\underline Y)=\lim A(Y_n)
=\lim \int_0^{T_n} L(q_n(s),\dot q_n(s),t+s)
ds
=h(q,t;q',t').
$$
As a consequence, we have
$$
\tilde h(P,t,-H(P,t));P',t',-H(P',t'))
\leq h(q,t;q',t').
$$
This ends the proof of the first part of the Proposition.

Let now $Y=(P,t,E)\in T_q^*M\times T^*\Tm$ and
 $Y'=(P',t',E')\in T_{q'}^*M\times T^*\Tm$ be
points such that
  $h(q,t;q',t') =\tilde h(Y;Y')$.
Let $q(s)$ be the projection on $M$ of the orbit
$\Phi_s(Y)$.
Using \ref{invariance} and \ref{Peierl}, we get 
$$
\tilde h(Y, Y')
=
\tilde h(\Phi_s(Y),Y')+\int _0^s
\lambda _{\Phi_{\sigma}(Y)}(V_G(\Phi_{\sigma}(Y))-G(\Phi_{\sigma}(Y))
d\sigma
$$
$$
\geq h(q(s),t+s;q',t')+\int_0^s L(q(\sigma),\dot q(\sigma),t+\sigma)dt \geq
h(q,t;q',t')= \tilde h(Y,Y').
$$
As a consequence, all the inequalities are equalities.
We obtain that the curve
$q(s)$ is minimizing in the expression
$$
h(q,t;q',t')=
\min\left(
h(q(s),t+s;q',t')+\int_0^s L(q(\sigma),\dot q(\sigma),t+\sigma)dt 
\right).
$$
Fathi has proved that $-P$ is then a super-differential
of the function $h(.,t;q',t')$ at $q$.
The properties at $(q',t')$ are treated in a similar way.
\qed
\subsection{}
\textsc{Corollary}
\begin{itshape}
If $H$ satisfies the hypotheses of 
\ref{HypothesesH}, then $m(H)\leq \tilde m(H)$.
\end{itshape}
\subsection{}
\textsc{Corollary}
\begin{itshape}
If $H$ satisfies the hypotheses of 
\ref{HypothesesH}, and if $m(H)=0$, then
$\tilde m(H)=0$, and we have 
 have $\tilde \mA_s(G)=\tilde \mA(H)$.
In addition,
we have 
$$
\tilde h(X_0,t_0,E_0;X_1,t_1,E_1)=
h(\pi(P_0),t_0;\pi(P_1),t_1)
$$
for each $(P_0,t_0,E_0)$ and $(P_1,t_1,E_1)$ in $\tilde \mA(H)$.
\vs
\end{itshape}

\proof
Let $(P,t,E)$ be a point of $T^*N$ and $q=\pi(P)$.
If $(P,t,E)\in \tilde \mA_s(G)$, then  
$$
\tilde h(P,t,E;P,t,E)=0,
$$
so that  $h(q,t;q,t)\leq 0$.
Since, on the other hand, we have $h(q,t;q,t)\geq m(H)=0$, we conclude that 
$h(q,t;q,t)=0$  hence 
$(q,t)\in \mA(H)$.
As a consequence, the function $h(q,t;.,t)$ 
is differentiable at $q$, see \ref{Aubry}, and  
$\big(\partial_3h(q,t;q,t),t
-H(\partial_3h(q,t;q,t),t)\big)\in \tilde \mA(H)$. 
Since $\tilde h(P,t,E;P,t,E)=h(q,t;q,t)$, the point
$P$ is a super-differential of $h(q,t;.,t)$ at
$q$, and we must have $P=\partial_3h(q,t;q,t)$.
Moreover, we have $G(P,t,E)=H(P,t)+E=0$, hence $(P,t,E)\in \tilde \mA(H)$.

Conversely, assume that $(P,t,E)\in \tilde \mA(H)$.
We then have $E=-H(P,t)$.
In addition,  $h(q,t;q,t)=0$,
the functions $h(q,t;.,t)$ and $h(.,t;q,t)$ are differentiable at $q$,
and we have $P=\partial_3h(q,t;q,t)=-\partial_1h(q,t;q,t)$.
Now let $X\in T_q^*M$
and $X'\in T_{q}^*M$ be such that 
$$
\tilde h(X,t,-H(X,t);X',t',-H(X',t'))=h(q,t;q,t).
$$
Then $-X$ is a super-differential at $q$ of $h(.,t;q,t)$,
and $X'$ is a super-differential at $q$ of $h(q,t;.,t)$.
It follows that $X=P=X'$.
Hence we have 
$\tilde h(P,t,E;P,t,E)=h(q,t;q,t)=0$.
This proves that $\tilde m(H)=0$, and that 
$(P,t,E)\in \tilde \mA_s(G)$.

 Finally, let $(P_0,t_0,E_0)\in T_{q_0}^*M\times T^*\Tm$ and 
$(P_1,t_1,E_1)\in T_{q_1}^*M\times T^*\Tm$ be two points of $\tilde \mA(H)$.
We have $E_0=-H(P_0,t_0)$ and $E_1=-H(P_1,t_1)$.
Furthermore,  the function
$h(q_0,t_0;.,t_1)$ is differentiable at $q_1$, with 
$\partial_3h(q_0,t_0;q_1,t_1)=P_1$,
and that the function $h(.,t_0;q_1,t_1)$ is differentiable at $q_0$,
with 
$\partial_1h(q_0,t_0;q_1,t_1)=-P_0$.
Since $-P_0$ and $P_1$ are then the only super-differentials
of  $h(.,t_0;q_1,t_1)$ and $h(q_0,t_0;.,t_1)$, we conclude that 
$\tilde h(P_0,t_0,E_0;P_1,t_1,E_1)=h(q_0,t_0;q_1,t_1)$.
\qed

\subsection{}
\textsc{Corollary}
\begin{itshape}
If $H$ satisfies the hypotheses of 
\ref{HypothesesH}, and if $m(H)=0$, then
 $\tilde \mM_s(G)=\tilde \mM(H)$.
\end{itshape}

\subsection{}
\textsc{Corollary}
\begin{itshape}
If $H$ satisfies the hypotheses of 
\ref{HypothesesH}, and if $m(H)=0$, then
 $\tilde \mN_s(G)=\tilde \mN(H)$.\vs
\end{itshape}

\proof
It is enough to prove that, if $(P_0,t_0,E_0)$
and $(P_1,t_1,E_1)$ belong to $\tilde \mA_s(G)$,
and $q_0=\pi(P_0), q_1=\pi(P_1)$,
then
$$\tilde \mI_s(P_0,t_0,E_0;
P_1,t_1,E_1)
=
\tilde \mI(q_0,t_0,q_1,t_1).
$$
Let $(P,t,E)$ be a point of $\tilde \mI_s(P_0,t_0,E_0;
P_1,t_1,E_1)$. 
We then have $G(P_0,t_0,E_0)=G(P,t,E)=0$
hence $E=-H(P,t)$.
Furthermore, the inequalities  
$$
h(q_0,t_0;q_1,t_1)=
\tilde h(P_0,t_0,-H(P_0,t_0);P_1,t_1,-H(P_1,t_1))
$$
$$
=\tilde h(P_0,t_0,-H(P_0,t_1);P,t,E)+
\tilde h(P,t,E; P_1,t_1,-H(P_1,t_1))
$$
$$
\geq h(q_0,t_0;q,t)+h(q,t;q_1,t_1)\geq h(q_0,t_0;q_1,t_1)
$$
are all equalities.
As a consequence, the point $(q,t)$ belongs 
to the set 
 $\mI(q_0,t_0;q_1,t_1)$,
and the differentials
$\partial_3h(q_0,t_0;q,t)$ and $\partial_1h(q,t;q_1,t_1)$ exist,
we have  $\partial_3h(q_0,t_0;q,t)=-\partial_1h(q,t;q_1,t_1)$,
and  the point 
$$(X,t,e)=(\partial_3h(q_0,t_0;q,t),t
,-H(\partial_3h(q_0,t_0;q,t),t))
$$ 
belongs
to $\tilde \mI(q_0,t_0;q_1,t_1)$, as follows from our definition
of the Ma\~n\'e set.
Since 
$$
\tilde h(P_0,t_0,-H(P_0,t_0);P,t,-H(P,t))=h(q_0,t_0;q,t),
$$ 
the point $P$ must
be a super-differential of $h(q_0,t_0;.,t)$ at $q$, hence
$P=X$.
We have proved that $(P,t,E)\in\tilde \mI(q_0,t_0;q_1,t_1)$.

Conversely, assume that $(P,t,E)\in \tilde \mI(q_0,t_0;q_1,t_1)$,
so that $E=-H(P,t)$.
Then  
$$
h(q_0,t_0;q,t)+h(q,t;q_1,t_1)=h(q_0,t_0;q_1,t_1)
$$
and 
$$
P=\partial_3h(q_0,t_0;q,t)=-\partial_1h(q,t;q_1,_1).
$$
In addition, since $(q_0,t_0)$ and $(q_1,t_1)$ belong to $\mA(H)$, 
the differential
$P_0=\partial_1h(q_0,t_0;q,t)$ exists for all $q$, 
and satisfies $(P_0,t_0,-H(P_0,t_0))\in \tilde \mA(H))$.
Similarly, setting
$P_1=\partial_3h(q,t;q_1,t_1)$, we have
$(P_1,t_1,-H(P_1,t_1))\in \tilde \mA(H))$.
We conclude that
$$
\tilde h(P_0,t_0,-H(P_0,t_0);P,t,E)=h(q_0,t_0;q,t)
$$
and 
$$
\tilde h(P,t,E;P_1,t_1,-H(P_1,t_1))= h(q,t;q_1,t_1).
$$
As a consequence, setting $E_0= -H(P_0,t_0)$ and $E_1=-H(P_1,t_1)$, we have
$$
\tilde h(P_0,t_0,E_0;P,t,E)
+\tilde h(P,t,E;P_1,t_1,E_1)
=
\tilde h(P_0,t_0,E_0;P_1,t_1,E_1).
$$
\qed

\begin{small}

\end{small}
\end{document}